\documentclass[12pt,a4paper]{article}
\usepackage{amsmath,amssymb,amsfonts}
\def\Bbb{\mathbb}

\title{\bf The number of Dedekind sums with equal fractional parts}

\author{Kurt Girstmair}
\date{}
\makeatletter
\let\@@maketitle=\maketitle
\def\maketitle{\def\thispagestyle##1{\relax}\@@maketitle}
\makeatother
\newtheorem{theorem}{Theorem}

\def\BE{\begin{equation}}
\def\EE{\end{equation}}
\def\BD{\begin{displaymath}}
\def\ED{\end{displaymath}}
\def\BA{\begin{array}}
\def\EA{\end{array}}
\def\BEA{\begin{eqnarray*}}
\def\EEA{\end{eqnarray*}}
\def\BI{\bibitem}

\def\Z{\Bbb Z}

\def\phi{\varphi}
\def\EPS{\varepsilon}

\def\CMOD#1#2#3{#1 \equiv #2 \: \mbox{mod}\: #3}

\def\MB{\mbox}
\def\LD{\ldots}
\def\OV{\overline}

\def\DIV{\,|\,}
\def\NDIV{\, \nmid \,}

\def\MN{\medskip\noindent}

\def\DED{Dedekind }

\begin{document}
\maketitle

\begin{abstract}
\noindent
In the paper \cite{Gi} it was shown that the \DED sums $12s(m,n)$ and $12s(x,n)$, $1\le m,x\le n$,
$(m,n)=(x,n)=1$,
are equal mod $\Z$ if, and only if,
$(x-m)(xm-1)\equiv 0$ mod $n$. Here we determine the cardinality of numbers $x$ in the above range that satisfy
this congruence for a given number $m$.
\end{abstract}

%%%%%%%%%%%%%%%%%%%%%%%%%%%%%%%%%%%%%%
\section*{1. Introduction and results}
%%%%%%%%%%%%%%%%%%%%%%%%%%%%%%%%%%%%%%

Let $n$ be a positive integer and $m\in \Z$, $(m,n)=1$. The classical \DED sum $s(m,n)$ is defined by
\BD
   s(m,n)=\sum_{k=1}^n ((k/n))((mk/n))
\ED
where $((\LD))$ is the usual sawtooth function (see, for instance, \cite{RaGr}, p. 1). In the present setting it is more
natural to work with
\BD
S(m,n)=12s(m,n)
\ED instead. Observe that $S(m+n,n)=S(m,n)$, hence it suffices to consider
arguments $m$ in the range $1\le m\le n$.
In the paper \cite{Gi} the following theorem was shown:

\begin{theorem} % Theorem 1 %%%%%%%%%%%%%%%%%%%%%%%%%%%%%%%%%%%%%%%%%%%%%%%%%%%%%
\label{t1}
Let $m, x$ be integers that are relatively prime to $n$. Then $S(x,n)\in S(m,n)+\Z$ if, and only if,
\BE\
\label{2}
(x-m)(xm-1)\equiv 0 \MB{ mod }n.
\EE
\end{theorem} %%%%%%%%%%%%%%%%%%%%%%%%%%%%%%%%%%%%%%%%%%%%%%%%%%%%%%%%%%%%%%%%%%%%%

\MN
For the history of this result, see \cite{Gi}. 
By means of Theorem \ref{t1} we can easily decide whether a \DED sum $S(x,n)$ has the same fractional part as
a given \DED sum $S(m,n)$. Therefore, it seems natural to count the respective arguments $x$, i.e., to
determine the number
\BD
\label{4}
  L(m,n)=|\{x\in \Z: 1\le x\le n, (x,n)=1, x \MB{ satisfies } (\ref{2})\}|
\ED
for a given integer $m$ with $(m,n)=1$. Let
\BD
  n=\prod_{p|n}p^{k_p}
\ED
be the prime decomposition of $n$. Obviously, the Chinese remainder theorem
implies
\BD
 L(m,n)=\prod_{p|n}L(m,p^{k_p}).
\ED
Hence it suffices to determine the numbers $L(m,p^k)$ for all primes $p$ not dividing $m$ and all positive integers $k$.
In this paper we show

\begin{theorem} % Theorem 2 %%%%%%%%%%%%%%%%%%%%%%%%%%%%%%%%%%%%%%%%%%%%%%%%%%%%%
\label{t2}
Let $p$ be a prime $\ge 3$, $k$ a positive integer, $m$ an integer not divisible by $p$, $1\le m\le p^k$. Then
\BD
\label{6}
  L(m,p^k)=\left\{\begin{array}{ll}
                    2 & \MB{ if } m\not\equiv \pm 1 \MB{ mod } p;\\
                    \rule{0mm}{5mm}
                    p^{\lfloor k/2\rfloor} & \MB{ if } m\equiv \pm 1 \MB{ mod } p^{\lceil k/2\rceil}; \\
                    \rule{0mm}{5mm}
                    2p^j & \MB { if } m=\pm 1 + p^jr, 1\le j<k/2, p\NDIV r.
                  \end{array}\right.
\ED
\end{theorem} %%%%%%%%%%%%%%%%%%%%%%%%%%%%%%%%%%%%%%%%%%%%%%%%%%%%%%%%%%%%%%%%%%%%%%%

The case $p=2$ is settled by
\begin{theorem} % Theorem 3 %%%%%%%%%%%%%%%%%%%%%%%%%%%%%%%%%%%%%%%%%%%%%%%%%%%%%
\label{t3}
Let $k$ be a positive integer and $m$ an odd integer, $1\le m\le 2^k$. Then
\BD
\label{7}
  L(m,2^k)=\left\{\begin{array}{ll}
                    2^{\lfloor k/2\rfloor} & \MB{ if } m\equiv \pm 1 \MB{ mod } 2^{\lceil k/2\rceil}; \\
                    \rule{0mm}{5mm}
                    2^{j+2} & \MB { if } m=\pm 1 + 2^jr, 2\le j<k/2-1, 2\NDIV r;\\
                    \rule{0mm}{5mm}
                    2^{j+1} & \MB{ if } m=\pm 1+ 2^{j}r, j=k/2-1\ge 2, 2\NDIV r;\\
                    \rule{0mm}{5mm}
                    2^{j} & \MB{ if } m=\pm 1+ 2^{j}r,  j=(k-1)/2\ge 2, 2\NDIV r.
                  \end{array}\right.
\ED
\end{theorem} %%%%%%%%%%%%%%%%%%%%%%%%%%%%%%%%%%%%%%%%%%%%%%%%%%%%%%%%%%%%%%%%%%%%%%%

Note that the case $j=1$ need not be considered in Theorem \ref{t3}. Indeed, suppose $m=1+2r$
for an odd number $r$. Then $\CMOD m34$, i.e., $m=-1+2^jr'$, $j\ge 2$, $r'\in\Z$, $r'$ odd.

\MN
{\em Example}: Let $n=1728=2^63^3$ and $m=7$. Since $\CMOD m{-1}4$ and $\CMOD m13$, our theorems give
$L(7,1728)=L(7,64)\cdot L(7,27)=8\cdot 6=48$. So we have 48 numbers $x$, $1\le x\le 1728$, $(x,6)=1$, such that
$S(x,1728)$ has the same fractional part $q$ as $S(7,1728)$, namely, $q=127/864$. These 48 \DED sums take the
positive values $q$, $3+q$, $8+q$, $11+q$, $16+q$, $24+q$, $31+q$, $56+q$, $59+q$, $248+q$,
and the negative values $-1+q$, $-5+q$, $-8+q$, $-9+q$, $-21+q$, $-32+q$.

\section*{2. Proof of Theorem \ref{t2}}

Let $p$ be a prime $\ge 3$ and $m$ an integer, $p\NDIV m$. For a given positive integer $k$, we simply write
\BD
  L=L(m,p^k).
\ED
Suppose that the integer $x$ lies in the range $1\le x\le p^k$. If $x$ satisfies the congruence (\ref{2}), then
we either have $\CMOD xmp$ or $\CMOD {xm}1p$, so $p\NDIV x$ in both cases. Hence
\BD
 L=|\{x\in \Z: 1\le x<p^k, x \MB{ satisfies } (\ref{2})\}|.
\ED
The case $m\not\equiv \pm 1$ mod $p$ is easy. Indeed, either $\CMOD xmp$ or $\CMOD{xm}1p$, and since
$p\ge 3$ each of these two
cases excludes the other one. So we have either $\CMOD xm{p^k}$ or $\CMOD{xm}1{p^k}$, and each of these two congruences
has exactly one solution $x$ with $1\le x\le p^k$. Therefore, we henceforth suppose
\BD
\label{8}
  \CMOD m{\EPS}p
\ED
with  $\EPS\in\{\pm 1\}$.

We consider the case $\CMOD m{\EPS}{p^{\lceil k/2\rceil}}$ first.  We show that $x$ satisfies (\ref{2}) if, and only if,
\BE
\label{9}
\CMOD x{\EPS}p^{\lceil k/2\rceil}.
\EE
If this is true, all possible solutions $x$ are given by
\BD
  x=\EPS+p^{\lceil k/2\rceil}s,
\ED
with
\BD
\left\{\begin{array}{ll}
         0\le s<p^{k-\lceil k/2\rceil}, & \MB{ if } \EPS=1; \\
         \rule{0mm}{5mm}
         0<s\le p^{k-\lceil k/2\rceil}, & \MB{ if } \EPS=-1.
       \end{array}\right.
\ED
Since $k-\lceil k/2\rceil=\lfloor k/2\rfloor$, there are $p^{\lfloor k/2\rfloor}$ suitable numbers $s$, and so
$L=p^{\lfloor k/2\rfloor}$.

In what follows we frequently consider the {\em $p$-exponent} $v_p(z)$ of an integer $z$,
namely,
\BD
  v_p(z)=\max\{j: p^j\DIV z\}.
\ED

In order to show the equivalence of (\ref{2}) and (\ref{9}) in the present case, we
first suppose $\CMOD x{\EPS}{p^{\lceil k/2\rceil}}$.
Since $\CMOD xm{p^{\lceil k/2\rceil}}$, $v_p(x-m)$  satisfies
\BD
  v_p(x-m)\ge \lceil k/2\rceil.
\ED
Further, if we write $m=\EPS +p^{\lceil k/2\rceil}r$, $x=\EPS+p^{\lceil k/2\rceil}s$ with $r,s\in \Z$,
we obtain
\BD
 xm-1=\EPS(r+s)p^{\lceil k/2\rceil}+p^{2\lceil k/2\rceil}rs
\ED
and $v_p(xm-1)\ge \lceil k/2\rceil$. Altogether $v_p((x-m)(xm-1))\ge 2\lceil k/2\rceil\ge k$,
so $x$ satisfies (\ref{2}).

Conversely, suppose $x\not\equiv\EPS$ mod $p^{\lceil k/2\rceil}$, so $v_p(x-\EPS)<k/2$. Since
$v_p(m-\EPS)\ge k/2$, we obtain
\BD
v_p(x-m)=v_p((x-\EPS)-(m-\EPS))=v_p(x-\EPS)<\lceil k/2\rceil.
\ED
Because $m=\EPS+p^{\lceil k/2\rceil}r$ and
$x=\EPS+p^ls$, $l<\lceil k/2\rceil$, $p\NDIV s$, we have
\BD
  xm-1=\EPS s p^l+(\EPS r+rsp^l)p^{\lceil k/2\rceil},
\ED
which shows $v_p(xm-1)=l<\lceil k/2\rceil$. Accordingly, $v_p((x-m)(xm-1))<2k/2=k$, and $x$ does not satisfy
(\ref{2}).

\medskip
\noindent{\em Remark.} In the case $\CMOD m{\EPS}p^{\lceil k/2\rceil}$, we have not used the condition $p\ge 3$,
so the assertion
\BD
   L(m,2^k)=2^{\lceil k/2\rceil}
\ED
is also true. Therefore, this case need not be considered in the proof of Theorem \ref{t3} below.

\medskip
\noindent
In the remainder of this section we suppose
\BD
  m=\EPS+p^jr, \enspace 1\le j<k/2, p\NDIV r.
\ED
If $x$ is a solution of (\ref{2}), it satisfies $\CMOD x{\EPS}p$.
%If $x=\EPS$, we have $v_p(x-m)=j=v_p(xm-1)$, so $x$ cannot satisfy (\ref{2}).
Hence we may write
\BD
\label{10}
 x=\EPS+p^ls, \enspace 1\le l\le k, p\NDIV s, 0< s\le p^{k-l}.
\ED
We first exclude the case $l\ne j$. Indeed, if this is true, we obtain $v_p(x-m)=\min\{j,l\}=v_p(xm-1)$.
However, $\min\{j,l\}\le j<k/2$, which means that $x$ is not a solution of (\ref{2}).
Accordingly, we may assume
\BE
\label{12}
x=\EPS+p^js, \enspace p\NDIV s, 0< s\le p^{k-j}
\EE
for a solution $x$ of (\ref{2}) in the range $1\le x \le n$.
We have
\BE
\label{14}
 x-m=p^j(s-r)\enspace \MB{ and }\enspace xm-1=p^j(\EPS(r+s)+p^jrs).
\EE

{\em Case 1:} $\CMOD{r-s}0p$. Since $p\ge 3$, $r+s\equiv 2r\not\equiv 0$ mod $p$. So (\ref{14}) yields
$v_p(xm-1)=j$. By (\ref{14}), $x$ is a solution of (\ref{2}) if, and only if, $v_p(s-r)\ge k-2j$.
Observe that $k-2j>0$,
since $j<k/2$. Accordingly, the parameter $s$ of a solution $x$ can be written as
\BE
\label{16}
 s=r+p^{k-2j}u,\enspace u\in\Z.
\EE
By (\ref{12}),
\BD
 0<r+p^{k-2j}u\le p^{k-j}, \MB{ i.e.}, \frac{-r}{p^{k-2j}}<u\le \frac{-r}{p^{k-2j}}+p^j.
\ED
Therefore, the integer $u$ lies in the interval
\BD
 \left]\frac{-r}{p^{k-2j}}, \frac{-r}{p^{k-2j}}+p^j\right]
\ED
of length $p^j$. This interval contains exactly $p^j$ integers $u$. By (\ref{16}), a number $s$ belonging to
such an integer $u$ is not divisible by $p$. Further, it satisfies $0<s\le p^{k-j}$ and, thus,
actually gives rise to a solution $x$
of (\ref{2}) in the range $1\le x\le n$. In other words, Case 1 yields $p^j$ solutions of the desired form.

{\em Case 2:} $r-s\not\equiv 0$ mod $p$. If $x$ is a solution of (\ref{2}), we necessarily have $\CMOD{r+s}0p$,
for otherwise $v_p(x-m)=j=v_p(xm-1)$, by (\ref{14}). But then $x$ cannot satisfy (\ref{2}) since $j<k/2$.
Observing $v_p(xm-1)=j+v_p(\EPS(r+s)+p^jrs)$, the condition
\BE
\label{18}
  v_p(\EPS(r+s)+p^jrs)\ge k-2j
\EE
is necessary and sufficient for $x$ being a solution of (\ref{2}).

{\em Subcase} (a): $k-2j\le j$, i.e., $j\ge k/3$.
Here $v_p(\EPS(r+s)+p^jrs)\ge k-2j$ if, and only if $v_p(r+s)\ge k-2j$.
So the number s has the form $s=-r+p^{k-2j}u$, $u\in\Z$. For any such integer $u$ we have $p\NDIV s$, because $p\NDIV r$,
and $k-2j>0$. By (\ref{12}), $s$ must be in the range $0<s\le p^{k-2j}$, which means
\BD
 \frac r{p^{k-2j}}<u\le\frac r{p^{k-2j}}+p^j.
\ED
As in case 1, there are $p^j$ integers $u$ with this property, all of which give rise to appropriate
solutions of (\ref{2}).

{\em Subcase} (b): $k-2j>j$, i.e., $j<k/3$. Here it turns out that $x$ is a solution of (\ref{2}) only if $v_p(r+s)=j$.
Indeed, $v_p(r+s)>j$ implies $v_p(\EPS(r+s)+p^jrs)=j<k-2j$, and $v_p(r+s)<j$ gives $v_p(\EPS(r+s)+p^jrs)<j<k-2j$.
Hence we may write $s=-r+p^ju$ for an integer $u$, $p\NDIV u$. Now
\BD
 \EPS(r+s)+p^jrs=p^j((\EPS+p^jr)u-r^2).
\ED
For a solution $x$ of the desired form it is necessary and sufficient that
\BE
\label{20}
v_p((\EPS+p^jr)u-r^2)\ge k-3j
\EE
(observe $k-3j\ge 1$). Let $\rho\in\Z$ be a representative of the residue class
$\OV r^2\OV{\EPS+p^jr}^{\:-1}\in \Z/p^{k-3j}\Z$.
Since $p\NDIV r$ we have $p\NDIV \rho$.  Moreover, $\CMOD u{\rho}{p^{k-3j}}$
is equivalent to (\ref{20}). If we write $u=\rho+p^{k-3j}v$, $v\in\Z$, condition (\ref{12}) gives
\BD
 \frac{r-p^j\rho}{p^{k-2j}}<v\le \frac{r-p^j\rho}{p^{k-2j}}+p^j.
\ED
There are exactly $p^j$ integers $v$ in this range, all of which give rise to appropriate solutions $x$ of (\ref{2}).

Altogether, the cases 1 and 2 yield $L=2p^j$.

\section*{Proof of Theorem \ref{t3}}

Let $m$ be an odd integer and $k\ge 1$.
Here we write $L=L(m,2^k)$. Let $\EPS\in\{\pm 1\}$. The case $\CMOD m{\EPS\,}{2^{\lceil k/2\rceil}}$ has
been settled in the foregoing section. Thus we may assume
\BD
 m=\EPS+2^jr, \enspace 2\le j<k/2, 2\NDIV r.
\ED
Since a solution $x$ of (\ref{2}) is odd, we write
\BE
\label{21}
  x=\EPS+2^ls, \enspace 1\le l\le k, 2\NDIV s, 0<s\le 2^{k-l}.
\EE
As in the proof of Theorem \ref{t2}, one can exclude the case $l\ne j$. Therefore, (\ref{21}) comes down to
\BE
\label{22}
  x=\EPS+2^js, \enspace 2\NDIV s, 0<s\le 2^{k-j}.
\EE
We obtain
\BE
\label{23}
  x-m=2^j(s-r), xm-1=2^j(\EPS(r+s)+2^jrs).
\EE
Because both numbers $r$ and $s$ are odd, we write
\BE
\label{23.1}
  s=r+2^tu, \enspace t\ge 1, u \MB{ odd }.
\EE
This gives
\BE
\label{24}
 s-r= 2^tu,\enspace \EPS(r+s)+2^jrs=2(\EPS r+2^{t-1}\EPS u+2^{j-1}r^2+2^{j+t-1}ru).
\EE

{\em Case 1:} $t\ge 2$. Since $j\ge 2$, (\ref{24}) implies $v_2(\EPS(r+s)+2^jrs)=1$. By (\ref{23}) and (\ref{24}),
$x$ is a solution
of (\ref{2}) if, and only if, $j+t+j+1\ge k$, i.e., $t\ge k-2j-1$. Observe $k-2j-1\ge 0$ because  $j<k/2$.
In view of (\ref{23.1}), we may write, for a solution $x$ of (\ref{2}),
$x=\EPS+2^js$, with $s=r+2^{k-2j-1}v$, $v\in\Z$. The condition $t\ge 2$ requires $4\DIV v$, if $k-2j-1=0$,
i.e., $j=(k-1)/2$; and $2\DIV v$, if $k-2j-1=1$, i. e., $j=k/2-1$. Now $x$ is in the appropriate range if,
and only if $0<s\le 2^{k-j}$, which means
\BD
  \frac{-r}{2^{k-2j-1}}<v\le \frac{-r}{2^{k-2j-1}}+ 2^{j+1}.
\ED
If $j<k/2-1$, there are $2^{j+1}$ integers $v$ in this range. If $j=k/2-1$, we have to count the number of even
integers $v$ in this range, which is $2^j$. In the case $j=(k-1)/2$, we count the number of integers
$v$ in this range with $4\DIV v$, which is $2^{j-1}$. Altogether, we have found $2^{j+1}$ solutions $x$
if $j<k/2-1$; $2^j$ solutions, if $j=k/2-1$; and $2^{j-1}$ solutions, if $j=(k-1)/2$.

{\em Case 2:} $t=1$. This means that $s=r+2u$ with $u$ odd. By (\ref{23}), $v_2(x-m)=j+1$. Further,
\BD
 \EPS(r+s)+2^jrs=2(\EPS r+\EPS u +2^{j-1}r^2+2^jru).
\ED
In view of (\ref{24}), a necessary and sufficient condition for $x$ being a solution of (\ref{2}) is
\BE
\label{26}
v_2(\EPS r+2^{j-1}r^2+u(\EPS +2^jr))\ge k-2(j+1).
\EE
Here we observe that $k-(2j+1)\ge 0$ if, and only if, $j\le k/2-1$.
We treat the case $j<k/2-1$ first. Since $k-2(j+1)>0$, (\ref{26}) is equivalent to
\BE
\label{28}
  u(\EPS +2^jr)\equiv -(\EPS r +2^{j-1}r^2) \MB{ mod } 2^{k-2(j+1)}.
\EE
Since $j\ge 2$, the right hand side of this congruence is odd, which means that $u$ is automatically an
odd number (as desired).
So $u$ has the form $u=\rho+2^{k-2(j+1)}v$ with some odd number $\rho$ coming from (\ref{28}) and $v\in\Z$.
Now $s=r+2\rho+2^{k-2j-1}v$ must satisfy $0<s\le 2^{k-j}$ (recall (\ref{22})).
This means
\BD
\frac{-r-2\rho}{2^{k-2j-1}}<v\le \frac{-r-2\rho}{2^{k-2j-1}}+2^{j+1}.
\ED
Since we have $2^{j+1}$ integers $v$ in this range, we obtain $2^{j+1}$ solutions $x$ of the desired form.
If $j\ge k/2-1$, the inequality (\ref{26}) is always fulfilled, and so we are left with the condition
\BD
  0<r+2u\le 2^{k-j}, \MB{ i.e., } -r/2<u\le -r/2+2^{k-j-1}.
\ED
In the case $j=k/2-1$, $k-j-1$ equals $j+1$, and there are $2^j$ odd numbers $u$ in this range. If $j=(k-1)/2$,
$k-j-1$ equals $j$,
and there are $2^j$ odd numbers $u$ of this kind. As in case 1, the number of solutions $x$ of the desired
form is $2^{j+1}$, if $j<k/2-1$; $2^j$, if $j=k/2-1$; and $2^{j-1}$, if $j=(k-1)/2$. On combining the cases 1 and 2,
we obtain
\BD
 L=\left\{\begin{array}{ll}
            2^{j+2} & \MB{ if } 2\le j<k/2-1; \\
            2^{j+1} & \MB{ if } j=k/2-1\ge 2; \\
            2^j & \MB{ if } j=(k-1)/2\ge 2.
          \end{array}\right.
\ED

%%%%%%%%%%%%%%%%%%%%%%%%%%%%%%%%%%%%%%%%%%%%%%%%%%%%%
%%%%%%%%%%%%%%%%%%%%%%%%%%%%%%%%%%%%%%%%%%%%%%%%%%%%%%%%%%%%%%%%%%%%%%%%%%

\vspace{0.5cm}
\noindent
Kurt Girstmair            \\
Institut f\"ur Mathematik \\
Universit\"at Innsbruck   \\
Technikerstr. 13/7        \\
A-6020 Innsbruck, Austria \\
Kurt.Girstmair@uibk.ac.at

\end{document}